\documentclass[11pt]{amsart}

\usepackage[margin=1in]{geometry}
\usepackage{amsmath,amssymb,amsthm}
\usepackage{mathtools}
\usepackage{graphicx}
\usepackage{enumitem}
\usepackage[normalem]{ulem}
\usepackage[dvipsnames]{xcolor}
\usepackage{tikz}
\usepackage[colorlinks=true,linkcolor=blue,citecolor=blue,urlcolor=blue]{hyperref}
\usepackage[capitalize,nameinlink]{cleveref}

\theoremstyle{plain}
\newtheorem{theorem}{Theorem}[section]

\theoremstyle{definition}

\theoremstyle{remark}

\usetikzlibrary{calc}
\definecolor{hexc0}{HTML}{E15759}
\definecolor{hexc1}{HTML}{F28E2B}
\definecolor{hexc2}{HTML}{EDC948}
\definecolor{hexc3}{HTML}{59A14F}
\definecolor{hexc4}{HTML}{76B7B2}
\definecolor{hexc5}{HTML}{4E79A7}
\definecolor{hexc6}{HTML}{B07AA1}
\definecolor{gR}{HTML}{DC2626}
\definecolor{gB}{HTML}{2563EB}
\definecolor{gG}{HTML}{16A34A}
\definecolor{gY}{HTML}{FACC15}
\definecolor{gW}{HTML}{000000}
\tikzset{
  g449edge/.style={draw=black!35, line width=0.15pt},
  g449vtx/.style={draw=none, line width=0.2pt, opacity=.9},
  p27edge/.style={draw=black!35, line width=0.1pt},
}

\newcommand{\R}{\mathbb{R}}

\newcommand{\Q}{\mathbb{Q}}
\newcommand{\C}{\mathbb{C}}

\title[From two heptagrams to a 449-vertex 5-chromatic unit-distance graph]{From two heptagrams to\\a 449-vertex 5-chromatic unit-distance graph}
\author{Marijn J. H. Heule \and Benjamin Przybocki}
\date{\today}

\begin{document}

\begin{abstract}
We construct a 5-chromatic unit-distance graph on 449 vertices. Previously, the smallest known 5-chromatic unit-distance graph had 509 vertices. Our graph is a subgraph of a graph obtained by a simple greedy construction applied to a 14-vertex seed graph consisting of two heptagrams.
\end{abstract}

\maketitle

\section{Introduction}

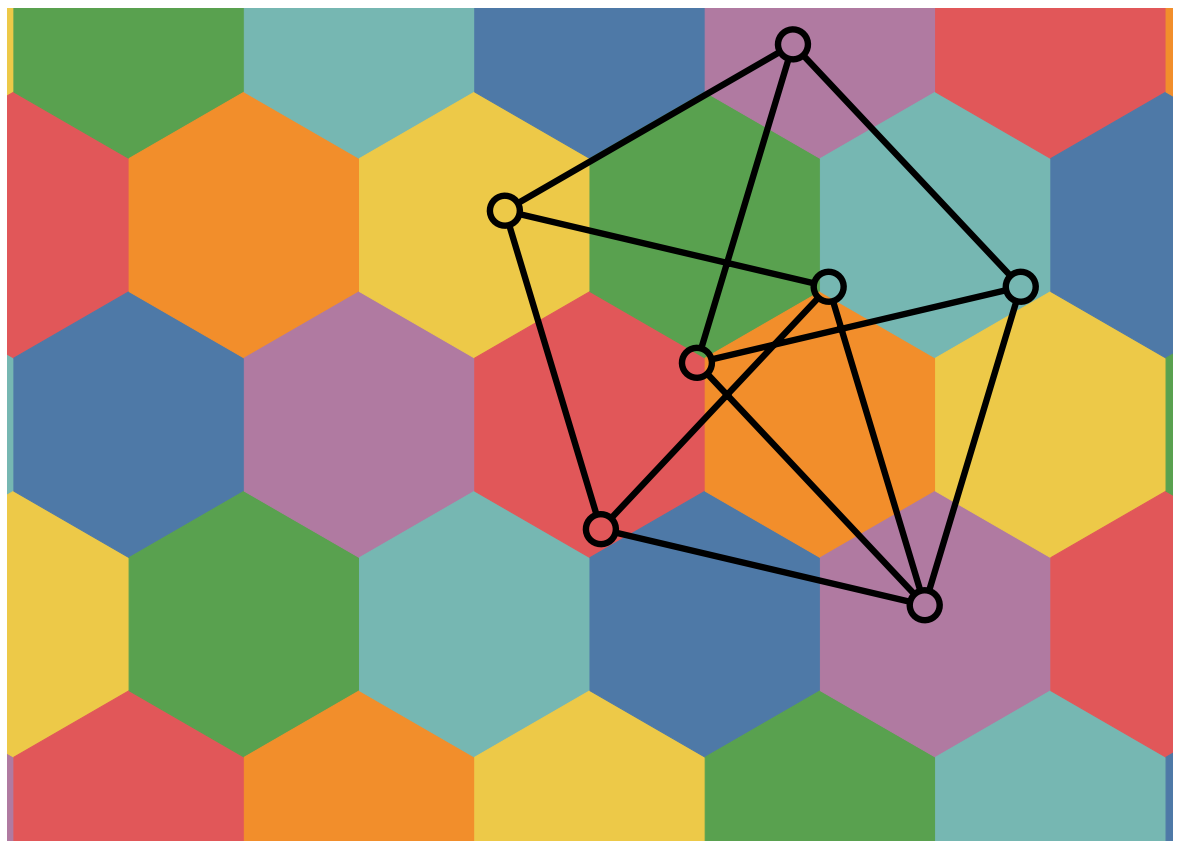
\begin{figure}[ht]
\centering
\begin{tikzpicture}[scale=2.2]
  \def\s{0.4}
  \clip (-1.75,-1.25) rectangle (1.75,1.25);
  \foreach \q in {-7,...,7} {
    \foreach \r in {-5,...,5} {
      \pgfmathtruncatemacro{\c}{mod(mod(\q+3*\r,7)+7,7)}
      \pgfmathsetmacro{\x}{\s*sqrt(3)*(\q+\r/2)}
      \pgfmathsetmacro{\y}{\s*1.5*\r}
      \filldraw[fill=hexc\c, draw=hexc\c, line width=0.2pt]
        (\x,\y) ++(90:\s) -- ++(-30:\s) -- ++(-90:\s) -- ++(-150:\s)
        -- ++(150:\s) -- ++(90:\s) -- cycle;
    }
  }
  \pgfmathsetmacro{\al}{asin(1/(2*sqrt(3)))}
  \begin{scope}[shift={(1.00636212,-0.54307033)}, rotate=30]
    \coordinate (a)  at (0,0);
    \coordinate (b1) at ({90-\al+30}:1);
    \coordinate (c1) at ({90-\al-30}:1);
    \coordinate (d1) at ({90-\al}:{sqrt(3)});
    \coordinate (b2) at ({90+\al+30}:1);
    \coordinate (c2) at ({90+\al-30}:1);
    \coordinate (d2) at ({90+\al}:{sqrt(3)});
    \begin{scope}[even odd rule]
    \clip (-5,-5) rectangle (5,5)
      (a) circle (0.045)
      (b1) circle (0.045) (c1) circle (0.045) (d1) circle (0.045)
      (b2) circle (0.045) (c2) circle (0.045) (d2) circle (0.045);
      \draw[line width=1.2pt, black, line cap=round, line join=round]
        (a) -- (b1) -- (d1) -- (c1) -- (a) (b1) -- (c1)
        (a) -- (b2) -- (d2) -- (c2) -- (a) (b2) -- (c2)
        (d1) -- (d2);
    \end{scope}
    \foreach \v in {a,b1,c1,d1,b2,c2,d2}
      \draw[draw=black, line width=1.2pt] (\v) circle (0.045);
  \end{scope}
\end{tikzpicture}
\caption{A $7$-coloring of the plane with no two points at distance $1$
receiving the same color, and the Moser spindle, a $4$-chromatic
unit-distance graph}
\label{fig:hex}
\end{figure}

How many colors are required to color $\mathbb{R}^2$ so that no pair of points at distance 1 from each other has the same color? Denote this number by $\chi(\mathbb{R}^2)$. This problem was posed by Edward Nelson in 1950 and has become known as the Hadwiger--Nelson problem (see \cite{soifer} for the history). A hexagonal tiling construction shows that $\chi(\mathbb{R}^2) \le 7$ (see \Cref{fig:hex}). Let a \emph{unit-distance graph} be a graph that can be embedded in the plane with each pair of adjacent vertices at distance 1 from each other. Then, the Moser spindle, also depicted in \Cref{fig:hex}, is a 4-chromatic unit-distance graph, which implies $\chi(\mathbb{R}^2) \ge 4$. In general, a compactness argument shows that $\chi(\mathbb{R}^2)$ is equal to the maximum chromatic number over all finite unit-distance graphs~\cite{compactness}, so determining $\chi(\mathbb{R}^2)$ reduces to a finitary combinatorial problem.

The bounds $4 \le \chi(\mathbb{R}^2) \le 7$ stood until 2018, when de Grey~\cite{degrey} found a 5-chromatic unit-distance graph on 1581 vertices, proving $\chi(\mathbb{R}^2) \ge 5$. Soon after, the collaborative Polymath16 project launched with the goal of finding smaller 5-chromatic unit-distance graphs~\cite{heule_553, heule_529, parts_509}. The culmination of this effort was a 5-chromatic unit-distance graph on 509 vertices found by Parts in 2020~\cite{parts_509}.

\begin{figure}[ht]
\centering
\begin{tikzpicture}[scale=2.0,
    vtx/.style={circle, fill=black, inner sep=1.3pt}]
  \pgfmathsetmacro{\rP}{1/(2*sin(180/7))}
  \pgfmathsetmacro{\rQ}{1/(2*sin(360/7))}
  \pgfmathsetmacro{\rR}{1/(2*sin(540/7))}
  \foreach \j in {0,...,6} {
    \coordinate (P\j) at ({360*\j/7+90}:\rP);
    \coordinate (Q\j) at ({360*\j/7+30}:\rQ);
    \coordinate (R\j) at ({360*\j/7+150}:\rR);
  }
  \foreach \j in {0,...,6} {
    \pgfmathtruncatemacro{\a}{mod(\j+1,7)}
    \pgfmathtruncatemacro{\b}{mod(\j+2,7)}
    \pgfmathtruncatemacro{\c}{mod(\j+3,7)}
    \draw[very thick, Plum] (P\j) -- (P\a);
    \draw[very thick, BurntOrange] (Q\j) -- (Q\b);
    \draw[very thick, OliveGreen] (R\j) -- (R\c);
    \draw[semithick] (P\j) -- (Q\j) -- (R\j) -- cycle;
  }
  \foreach \j in {0,...,6} \node[vtx] at (P\j) {};
  \foreach \j in {0,...,6} \node[vtx] at (Q\j) {};
  \foreach \j in {0,...,6} \node[vtx] at (R\j) {};
\end{tikzpicture}
\caption{Haugland's $H_{21}$, consisting of a regular heptagon and two regular heptagrams}
\label{fig:h21}
\end{figure}
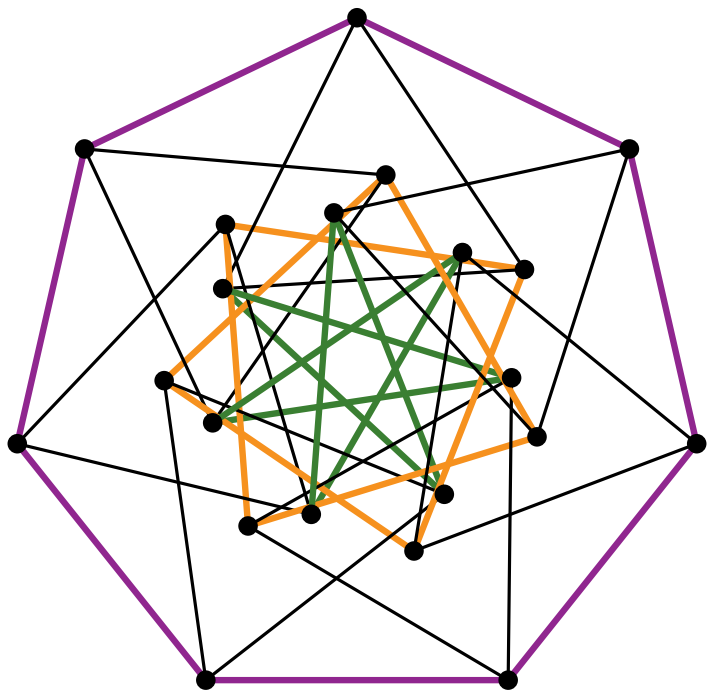

All of the above graphs use the Moser spindle as a building block. Voronov, Neopryatnaya, and Dergachev~\cite{no_moser} constructed a 5-chromatic unit-distance graph without the Moser spindle on 64,513 vertices, giving an alternative route to showing $\chi(\mathbb{R}^2) \ge 5$. Heule~\cite{heule_1441} improved on this by exhibiting a 5-chromatic unit-distance graph without the Moser spindle on 1441 vertices. Recently, Haugland~\cite{haugland} constructed a 5-chromatic unit-distance graph without the Moser spindle on 2131 vertices. While larger than Heule's graph, Haugland's construction is interesting in that, unlike earlier graphs, it isn't based on gluing together many copies of a small 4-chromatic unit-distance graph. Every edge in Haugland's graph is one of the unit vectors from the unit-distance graph $H_{21}$, shown in \Cref{fig:h21}. Since then, Gal\'an~\cite{galan} used a similar construction to obtain a 5-chromatic unit-distance graph without the Moser spindle on 852 vertices, and Shukla~\cite{shukla} found one on 557 vertices.

\begin{figure}[ht]
\centering
\resizebox{\textwidth}{!}{\input{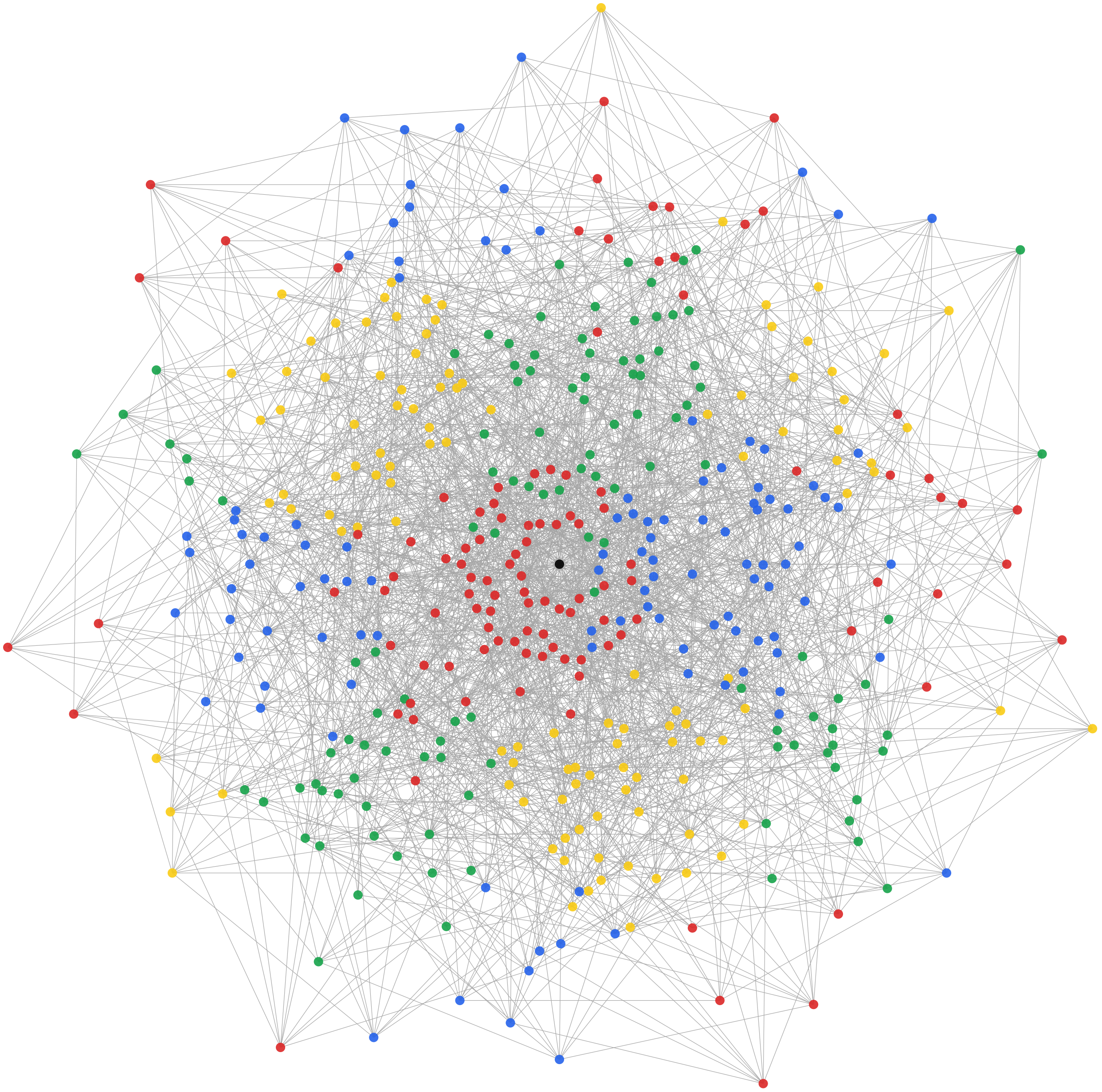}}
\caption{The graph $G_{449}$. The vertices are 5-colored in such a way that the center vertex is the only one colored black.}
\label{fig:g449}
\end{figure}

Building on Haugland's construction, we prove the following:
\begin{theorem} \label{thm:main}
    There is a 5-chromatic unit-distance graph on 449 vertices and 2198 edges.
\end{theorem}
We denote the graph witnessing \Cref{thm:main} by $G_{449}$, depicted in \Cref{fig:g449}. The graph has 7-fold rotational symmetry, and its coordinates lie in the field $\Q(\sin(2 \pi / 21))$; in particular, $G_{449}$ does not contain the Moser spindle as a subgraph.

We note that recently, OpenAI~\cite{OAI} proved that $\chi(\mathbb{R}^2) \ge 6$. Their proof is not constructive, and to the best of our knowledge, no finite 6-chromatic unit-distance graph is known.

\section{Constructing $G_{449}$}

We identify $\R^2$ with $\C$. With this identification, every point in $G_{449}$ lies in the cyclotomic field $\Q(e^{2\pi i/21})$. This implies that the real-valued coordinates lie in $\Q(\sin(2 \pi / 21))$.

At a high level, our construction proceeds in three steps. First, we start with a small \emph{seed graph}, a unit-distance graph that serves as the starting point of our construction. Then, we repeatedly add vertices to the seed graph according to a heuristic that tries to maximize edge count. Finally, once the graph is 5-chromatic, we find a small subgraph that is still 5-chromatic.

\subsection{Edge directions}

Each edge in $G_{449}$ belongs to a small set of unit vectors. Let $U$ be the set of $84$ unit vectors $b - a$, where $(a,b)$ ranges over the ordered pairs of adjacent vertices of $H_{21}$ (\Cref{fig:h21}), and let $\overline{U}$ be the set of their complex conjugates (i.e., their reflections across the $x$-axis). Then, every edge of $G_{449}$ is a unit vector from $D \coloneqq U \cup \overline{U}$. Note that since $H_{21}$ has 7-fold rotational symmetry, $D$ is invariant under rotation by $2\pi/7$ about the origin.

\subsection{Seed graph}

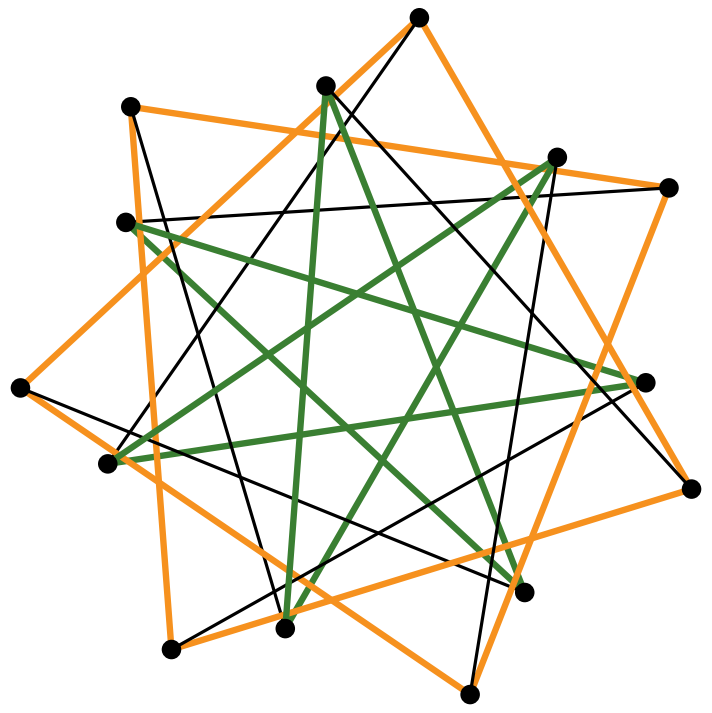
\begin{figure}[ht]
\centering
\begin{tikzpicture}[scale=3.6,
    vtx/.style={circle, fill=black, inner sep=1.3pt}]
  \pgfmathsetmacro{\rQ}{1/(2*sin(360/7))}
  \pgfmathsetmacro{\rR}{1/(2*sin(540/7))}
  \foreach \j in {0,...,6} {
    \coordinate (Q\j) at ({360*\j/7+30}:\rQ);
    \coordinate (R\j) at ({360*\j/7+150}:\rR);
  }
  \foreach \j in {0,...,6} {
    \pgfmathtruncatemacro{\b}{mod(\j+2,7)}
    \pgfmathtruncatemacro{\c}{mod(\j+3,7)}
    \draw[very thick, BurntOrange] (Q\j) -- (Q\b);
    \draw[very thick, OliveGreen] (R\j) -- (R\c);
    \draw[semithick] (Q\j) -- (R\j);
  }
  \foreach \j in {0,...,6} \node[vtx] at (Q\j) {};
  \foreach \j in {0,...,6} \node[vtx] at (R\j) {};
\end{tikzpicture}
\caption{The seed graph $C_{14}$}
\label{fig:c14}
\end{figure}

Our seed graph is the $14$-vertex induced subgraph $C_{14}$ of $H_{21}$ consisting of its two heptagrams (see \Cref{fig:c14}). With our coordinate system, its vertices are
\[
P_0 = \bigl\{ r_2\, e^{i(\pi/6 + 2\pi k/7)},\ r_3\, e^{i(5\pi/6 + 2\pi k/7)} \ :\ 0 \le k \le 6 \bigr\}, \qquad r_2 = \frac{1}{2\sin(2\pi/7)},\quad r_3 = \frac{1}{2\sin(3\pi/7)},
\]
which form a regular $\{7/2\}$ heptagram of radius $r_2$ and a regular $\{7/3\}$ heptagram of radius $r_3$. Together with the $7$ edges joining the two heptagrams, $C_{14}$ has $21$ edges.

\subsection{Growing procedure}

We grow the seed graph by repeating the following greedy operation. Given a point set $P_t$, let
\[
d_t(x) = \bigl|\{ p \in P_t : x - p \in D \}\bigr| \qquad \text{for } x \in (P_t + D) \setminus P_t,
\]
and simultaneously add all candidates $x$ attaining the maximum value of $d_t$:
\[
P_{t+1} = P_t \cup \bigl\{ x \in (P_t + D) \setminus P_t : d_t(x) = \textstyle\max_y d_t(y) \bigr\}.
\]
Since $P_0$ and $D$ are invariant under rotation by $2\pi/7$ about the origin, so is every $P_t$. Thus, each $P_t$ consists of 7-point orbits under the rotational symmetry, together with the origin once it is added (which happens at $t=4$).

\begin{figure}[ht]
\centering
\resizebox{\textwidth}{!}{\input{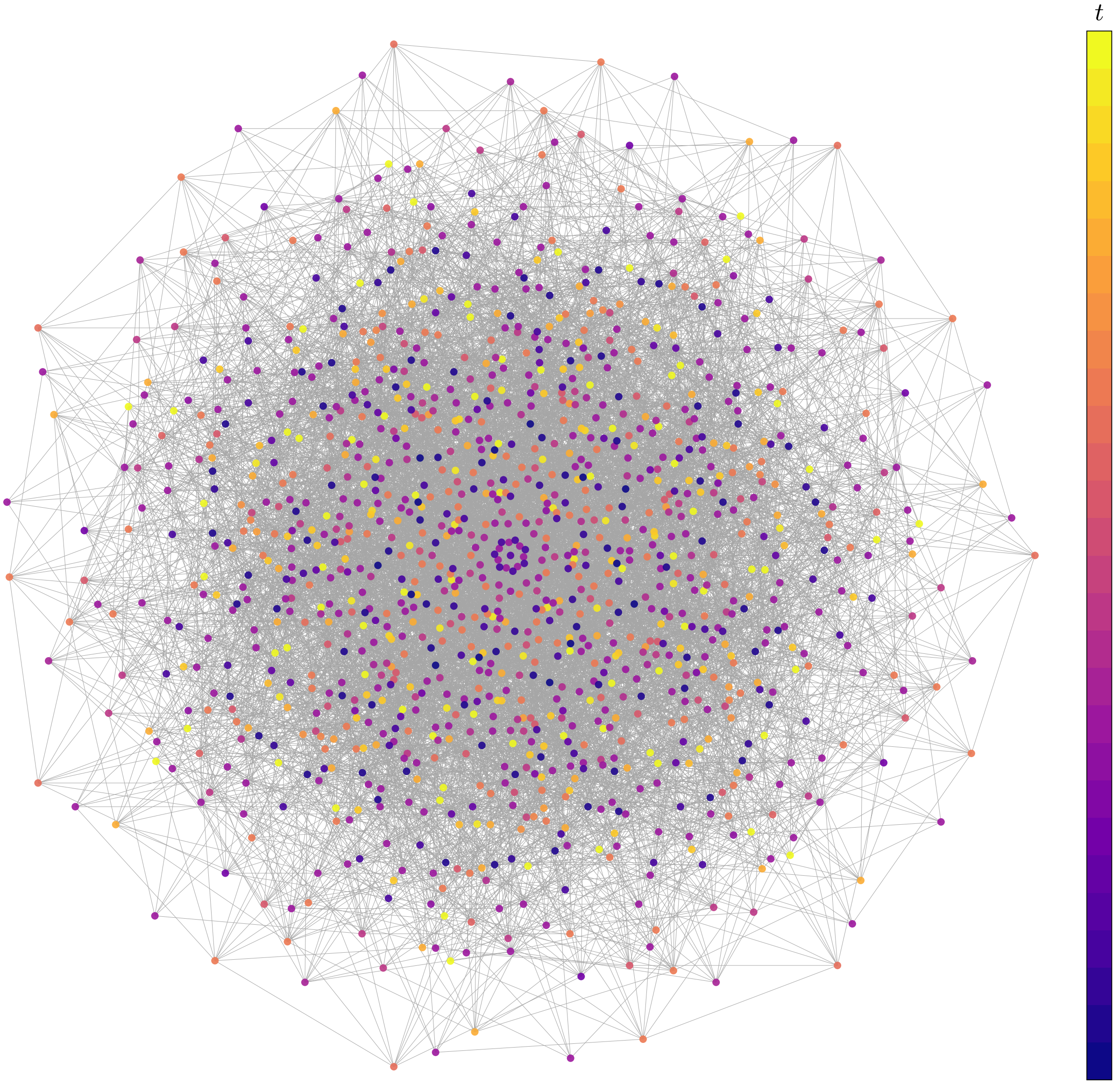}}
\caption{The graph $P_{27}$. Each vertex is colored according to the step $t$ at which it was added.}
\label{fig:p27}
\end{figure}

For each $t$, we identify $P_t$ with the unit-distance graph $(P_t, E_t)$, where $E_t$ is the set of pairs of points at distance 1 from each other. The graph $P_{17}$, which has 806 vertices, is the first 5-chromatic graph in the series. Nevertheless, we continue the growing procedure until $P_{27}$, which has 1233 vertices (see \Cref{fig:p27}). The additional points in $P_{27}$ give us more freedom in the reduction step. \Cref{tab:growth} shows some data for the graphs $P_t$. Note that $P_{27}$ without the origin vertex is 4-colorable, hence all $P_{t}$ with $t \leq 27$ are 5-colorable.

\begin{table}[ht]
\centering
\caption{The graphs $P_t$ of the growing procedure. For $t \ge 17$, the last two columns list the time required for \textsf{CaDiCaL} to prove that $P_t$ is not 4-colorable and the number of conflicts it produced, using the encoding of \Cref{sec:encoding} while fixing the colors in a triangle whose total degree is maximal.}
\label{tab:growth}
\begin{tabular}{rrrrrr}
\hline
$t$ & $|P_t|$ & $|E_t|$ & Triangles & Time (s) & Conflicts \\
\hline
0 & 14 & 21 & 0 & -- & -- \\
1 & 98 & 238 & 56 & -- & -- \\
2 & 112 & 308 & 70 & -- & -- \\
3 & 196 & 651 & 154 & -- & -- \\
4 & 197 & 658 & 154 & -- & -- \\
5 & 225 & 826 & 196 & -- & -- \\
6 & 239 & 896 & 210 & -- & -- \\
7 & 246 & 931 & 210 & -- & -- \\
8 & 260 & 1001 & 238 & -- & -- \\
9 & 582 & 2597 & 630 & -- & -- \\
10 & 610 & 2863 & 714 & -- & -- \\
11 & 652 & 3234 & 840 & -- & -- \\
12 & 715 & 3710 & 980 & -- & -- \\
13 & 722 & 3780 & 987 & -- & -- \\
14 & 736 & 3899 & 1029 & -- & -- \\
15 & 771 & 4179 & 1099 & -- & -- \\
16 & 785 & 4312 & 1141 & -- & -- \\
17 & 806 & 4480 & 1197 & 314 & $3{,}709{,}048$ \\
18 & 953 & 5558 & 1477 & 287 & $3{,}250{,}711$ \\
19 & 960 & 5642 & 1512 & 272 & $3{,}102{,}919$ \\
20 & 967 & 5712 & 1540 & 301 & $3{,}334{,}294$ \\
21 & 981 & 5845 & 1582 & 254 & $2{,}849{,}273$ \\
22 & 1051 & 6426 & 1750 & 282 & $2{,}947{,}257$ \\
23 & 1058 & 6496 & 1771 & 326 & $3{,}365{,}098$ \\
24 & 1121 & 7007 & 1904 & 311 & $3{,}128{,}470$ \\
25 & 1128 & 7084 & 1932 & 314 & $3{,}145{,}171$ \\
26 & 1156 & 7343 & 2016 & 269 & $2{,}816{,}153$ \\
27 & 1233 & 7973 & 2233 & 224 & $2{,}380{,}363$ \\
\hline
\end{tabular}
\end{table}

\subsection{Reduction}

The graph $G_{449}$ is the subgraph of $P_{27}$ induced by the origin and $64$ of its $176$ rotation orbits. We used a few heuristics to decide which orbits to remove.

First, we follow the approach of \cite{heule_529} by using \textsf{CaDiCaL}~\cite{cadical} to prove that $P_{27}$ (or some $P_t$ for $t \ge 27$) is not 4-colorable, trimming the proof with \textsf{drat-trim}~\cite{Wetzler2014, heule2016dratformatdrattrimchecker}, and computing an UNSAT core. This allows us to determine which vertices were used in \textsf{CaDiCaL}'s proof. We repeat this solve-and-trim procedure 20 times in sequence to obtain smaller UNSAT cores, and we rerun the entire procedure many times to obtain a variety of UNSAT cores. Then, we consider the subgraph obtained by keeping all orbits such that the average fraction of their vertices appearing in the UNSAT cores exceeds some threshold $p \in (0,1)$; intuitively, this selects the orbits that are most important to the graph's being 5-chromatic. We choose $p$ to be roughly as large as it can be so that the resulting subgraph is 5-chromatic.

Second, we remove orbits one at a time and check if the graph is still 5-chromatic. There are many choices of which orbits to remove, which results in several graphs.

At this point, we had several 5-chromatic unit-distance graphs on 456--484 vertices. We then considered hybrid graphs obtained by choosing sets of orbits that were present in some, but not all, of our graphs so that the total number of vertices was 456. This resulted in 171 distinct 5-chromatic unit-distance graphs on 456 vertices. For exactly one of these graphs, we were able to remove one orbit while preserving 5-chromaticity, resulting in $G_{449}$. We tried to shrink the graph even further by removing individual vertices rather than entire orbits, but every proper induced subgraph of $G_{449}$ is 4-colorable.

\subsection{Verification}

To prove \Cref{thm:main}, we must verify both that $G_{449}$ is a unit-distance graph and that it is 5-chromatic.

\subsubsection{Unit distances.}
To verify that $G_{449}$ is a unit-distance graph, we use exact arithmetic (i.e., no floating-point arithmetic). Specifically, let $\xi = e^{i\pi/21}$, a primitive $42$nd root of unity, so that $\Q(\xi) = \Q(e^{2\pi i/21})$. Every vertex can be written in the form $\frac{1}{49}\sum_{j=0}^{11} c_j \xi^j$ for some $c_0, \dots, c_{11} \in \mathbb{Z}$. Then, two vertices $u$ and $v$ are adjacent if and only if $d = u - v$ satisfies $d \overline{d} = 1$, which we test using integer arithmetic. By this means, we certify the 2198 edges of $G_{449}$.

\subsubsection{Encoding.}\label{sec:encoding}
We encode the 4-colorability of $G_{449}$ as a CNF formula as follows; our goal is to prove that the formula is unsatisfiable. For each vertex $v$ and color $c \in \{1,2,3,4\}$, we introduce a Boolean variable $x_{v,c}$ meaning that $v$ has color $c$. The formula consists of a clause $x_{v,1} \lor x_{v,2} \lor x_{v,3} \lor x_{v,4}$ for each vertex $v$ and a clause $\overline{x_{u,c}} \lor \overline{x_{v,c}}$ for each edge $\{u,v\}$ and color $c$. This formula is satisfiable if and only if $G_{449}$ is 4-colorable. Note that it is unnecessary to enforce that each vertex receives at most one color, because any satisfying assignment in which a vertex receives multiple colors can be transformed into one in which every vertex receives exactly one color.

To improve the performance of the encoding, we use clique hints~\cite{queens}, which are redundant constraints that help the solver. In a proper 4-coloring, the three vertices of a triangle receive three distinct colors. Thus, for every pair of colors, at least one of them appears in any given triangle. Accordingly, for each triangle $\{a,b,c\}$ of $G_{449}$ and each pair of colors $\{x,y\}$, we include the redundant clause $x_{a,x} \lor x_{a,y} \lor x_{b,x} \lor x_{b,y} \lor x_{c,x} \lor x_{c,y}$.

\subsubsection{Symmetry breaking and case splitting.}
The formula has two kinds of symmetry: permutations of the four colors and the 7-fold rotational symmetry of $G_{449}$. We break both using the neighborhood of the origin. We may assume that the origin has color $1$. Let its seven neighbors be $n_0, \dots, n_6$, numbered so that $n_{k+1}$ is the rotation of $n_k$ by $2\pi/7$ for each $k$. These vertices are pairwise non-adjacent and must receive colors in $\{2,3,4\}$. It is not necessary to consider all $3^7 = 2187$ possible colorings of $n_0, \dots, n_6$, because two such colorings are equivalent if one can be obtained from the other by a rotation and a permutation of the colors $\{2,3,4\}$. Some color appears at least three times among $n_0, \dots, n_6$; we may assume it is color $2$. Rotating if necessary, we may also assume that $n_0$ receives color $2$. Classifying the possible placements of the $2$s gives the six cases in \Cref{tab:cases}.

\begin{table}[ht]
\centering
\renewcommand{\ULthickness}{0.8pt}
\caption{The six cases of the symmetry breaking. In a pattern, an uncrossed digit is a fixed color, \texttt{\protect\sout{2}} means ``not color $2$'', and \texttt{.} is unconstrained. For each case, we list the time required for \textsf{CaDiCaL} to solve the case as well as the number of conflicts it produced.}
\label{tab:cases}
\begin{tabular}{clrr}
\hline
Case & Pattern $n_0 \cdots n_6$ & Time (s) & Conflicts \\
\hline
1 & \texttt{2222222} & $15$ & $322{,}711$ \\
2 & \texttt{2223...} & $3505$ & $27{,}961{,}976$ \\
3 & \texttt{2232..\sout{2}} & $2316$ & $19{,}696{,}966$ \\
4 & \texttt{223\sout{2}2\sout{22}} & $784$ & $8{,}051{,}763$ \\
5 & \texttt{223\sout{22}2\sout{2}} & $874$ & $8{,}932{,}088$ \\
6 & \texttt{232\sout{2}2\sout{22}} & $772$ & $8{,}340{,}369$ \\
\hline
Total & & $8266$ & $73{,}305{,}873$ \\
\hline
\end{tabular}
\end{table}

To reduce the overlap between cases, we process them in the order $1, 2, 3, 5, 4, 6$ and, in each case, exclude every coloring of $n_0, \dots, n_6$ that is equivalent to a coloring admitted by a previous case or earlier in the same case. Each redundant coloring is excluded by an additional clause. Below, we describe how we proved that the resulting case split is exhaustive.

\subsubsection{Solving.}
With our encoding, each case in \Cref{tab:cases} corresponds to a CNF formula; $G_{449}$ is 4-colorable if and only if at least one of these formulas is satisfiable. We solved the six formulas with \textsf{CaDiCaL}~3.0.1~\cite{cadical}, using the options \texttt{--forcephase --scorefactor=950 --sat}, on an Apple M2 Pro. Despite the name, the \texttt{--sat} option makes the solving faster for our formulas. All six formulas are unsatisfiable. \Cref{tab:cases} lists the solving time and number of conflicts for each case.

Interestingly, \textsf{CaDiCaL} proves that Parts' 509-vertex graph~\cite{parts_509} is not 4-colorable in only about $3$ seconds and $70{,}000$ conflicts using a similar encoding. Thus, the absence of a 4-coloring is considerably harder to establish for $G_{449}$ than for Parts' graph.

\subsubsection{Proof checking.}
The goal of proof checking is to certify that (a) the symmetry breaking is sound, (b) the case split is exhaustive, and (c) the CNF formulas are indeed unsatisfiable. We certify all three at once with a single proof of unsatisfiability of the CNF given by the encoding from \Cref{sec:encoding} (without the clique hints) in the linear substitution redundancy (LSR) format~\cite{codel_sr}. The proof consists of three parts. First, a proof in the substitution redundancy (SR) proof system~\cite{codel_sr} adds the clique hints, the unit clause fixing the origin to have color $1$, and clauses excluding each of the colorings of $n_0, \dots, n_6$ not covered by one of the six cases in \Cref{tab:cases}. We used \textsf{dsr-trim}~\cite{codel_sr} to convert this proof to LSR. Second, for each case, we used \textsf{CaDiCaL} to produce a proof of unsatisfiability in the LRAT format~\cite{lrat}, and we added the negation of the case's assumptions to every clause derived in this proof; thus, we obtain a proof whose conclusion is the negation of the case in question. Third, we used \textsf{CaDiCaL} to prove that the clauses derived in the previous two parts are jointly unsatisfiable, which is to say that the case split is exhaustive. We combined these steps into a single LSR proof and streamed it into \textsf{lsr-check}~\cite{codel_sr}, which verified it. Therefore, $G_{449}$ is not 4-colorable. On the other hand, \Cref{fig:g449} shows that $G_{449}$ is 5-colorable (unsurprisingly), so $G_{449}$ is 5-chromatic, establishing \Cref{thm:main}.

The graph $G_{449}$, together with scripts for reproducing these verifications, is available in a GitHub repository~\cite{g449repo}.

\subsection*{Acknowledgments}

This research was supported by the DARPA expMath program through the DARPA CMO contract number HR0011262E028 and by the
National Science Foundation (NSF) under grant DMS2434625.
Przybocki was additionally supported by the NSF Graduate Research Fellowship Program under Grant Nos. DGE-2140739 and DGE-2631988.

\bibliographystyle{amsplain}
\bibliography{references}

\end{document}